\documentclass[english]{amsart}

\usepackage{amsmath} 
\usepackage{amssymb} 
\usepackage{amscd}
\usepackage{tabularx}

\title{A modularity lifting theorem for weight two Hilbert modular forms}
\author{Toby Gee}
\email{toby.gee@ic.ac.uk}
\address{Department of Mathematics, Imperial College London}


\newcommand{\F}{\mathbb{F}} 

 \newcommand{\M}{\mathcal{M}}

 \newcommand{\Z}{\mathbb{Z}}
 
 \newcommand{\Q}{\mathbb{Q}}

 \newcommand{\Gal}{\operatorname{Gal}}
 \newcommand{\GL}{\operatorname{GL}}
\newcommand{\FI}{\operatorname{FI}}

\newcommand{\rhobar}{\overline{\rho}}

\newcommand{\ad}{\operatorname{ad}}

\newcommand{\End}{\operatorname{End}}

\newcommand{\Rep}{\operatorname{Rep}}

\newcommand{\Mod}{\operatorname{Mod}}
\newcommand{\diag}{\operatorname{diag}}
\newcommand{\bigs}{\mathfrak{S}}

\usepackage{amsthm}

\newtheorem*{thmn}{Theorem}

 \newtheorem{thm}{Theorem}[section]
 
 \newtheorem{lemma}[thm]{Lemma}
 \newtheorem{prop}[thm]{Proposition}
 \theoremstyle{definition}
 
 \theoremstyle{remark}
 
 \numberwithin{equation}{subsection}

\begin{document}
\maketitle
\begin{abstract}We prove a modularity lifting theorem for potentially
Barostti-Tate representations over totally real fields, generalising recent
results of Kisin. Unfortunately, there was an error in the original version of this paper, meaning that we can only obtain a slightly weaker result in the case where the representations are potentially ordinary; an erratum has been added explaining this error.

\end{abstract}
\section{Introduction}In \cite{kis04} Mark Kisin introduced a number of new
techniques for proving modularity lifting theorems, and was able to prove a
very general lifting theorem for potentially Barsotti-Tate representations
over $\Q$. In \cite{kis05} this was generalised to the case of $p$-adic
representations of the absolute Galois group of a totally real field in
which $p$ splits completely. In this note, we further generalise this result
to:

\begin{thmn}Let $p>2$, let $F$ be a totally real field in which $p$ is unramified, and
let $E$ be a finite extension of $\Q_p$ with ring of integers $\mathcal{O}$.
Let $\rho:G_F\to\GL_2(\mathcal{O})$ be a continuous representation
unramified outside of a finite set of primes, with determinant a finite
order character times the $p$-adic cyclotomic character. Suppose that
\begin{enumerate}\item $\rho$ is potentially Barsotti-Tate at each $v|p$.
\item $\overline{\rho}$ is modular. \item $\overline{\rho}|_{F(\zeta_p)}$ is
absolutely irreducible.\end{enumerate}Then $\rho$ is modular.\end{thmn}

We emphasise that the techniques we use are entirely those of Kisin. Our
only new contributions are some minor technical improvements; specifically,
we are able to prove a more general connectedness result than Kisin for
certain local deformation rings, and we replace an appeal to a result of
Raynaud by a computation with Breuil modules with descent data.

The motivation for studying this problem was the work reported on in
\cite{gee053}, where we apply the main theorem of this paper to the
conjectures of \cite{bdj}. In these applications it is crucial to have a
lifting theorem valid for $F$ in which $p$ is unramified, rather than just
totally split.

\section{Connected components}\label{ctd}Firstly, we recall some definitions and theorems from
\cite{kis04}. We make no attempt to put these results in context, and the
interested reader should consult section 1 of \cite{kis04} for a more
balanced perspective on this material.

Let $p>2$ be prime. Let $k$ be a finite extension of $\F_p$ of cardinality
$q=p^r$, and let $W=W(k)$, $K_0=W(k)[1/p]$. Let $K$ be a totally ramified
extension of $K_0$ of degree $e$. We let $\mathfrak{S}=W[[u]]$, equipped
with a Frobenius map $\phi$ given by $u\mapsto u^p$, and the natural
Frobenius on $W$. Fix an algebraic closure $\overline{K}$ of $K$, and fix a
uniformiser $\pi$ of $K$. Let $E(u)$ denote the minimal polynomial of $\pi$
over $K_0$.

Let $'(\Mod/\bigs)$ denote the category of $\bigs$-modules $\mathfrak{M}$
equipped with a $\phi$-semilinear map $\phi:\mathfrak{M}\to\mathfrak{M}$
such that the cokernel of $\phi^*(\mathfrak{M})\to\mathfrak{M}$ is killed by
$E(u)$. For any $\Z_p$-algebra $A$, set $\bigs_A=\bigs\otimes_{\Z_p}A$.
Denote by $'(\Mod/\bigs)_A$ the category of pairs $(\mathfrak{M},\iota)$
where $\mathfrak{M}$ is in $'(\Mod/\bigs)$ and
$\iota:A\to\End(\mathfrak{M})$ is a map of $\Z_p$-algebras.

We let $(\Mod \FI/\bigs)_A$ denote the full subcategory of $'(\Mod/\bigs)_A$
consisting of objects $\mathfrak{M}$ such that $\mathfrak{M}$ is a
projective $\bigs_A$-module of finite rank.

Choose elements $\pi_n\in\overline{K}$ ($n\geq 0$) so that $\pi_0=\pi$ and
$\pi_{n+1}^p=\pi_n$. Let $K_{\infty}=\bigcup_{n\geq 1}K(\pi_n)$. Let
$\mathcal{O}_\mathcal{E}$ be the $p$-adic completion of $\bigs[1/u]$. Let
$\Rep_{\Z_p}(G_{K_\infty})$ denote the category of continuous
representations of $G_{K_\infty}$ on finite $\Z_p$-algebras. Let
$\Phi\operatorname{M}_{\mathcal{O}_\mathcal{E}}$ denote the category of
finite $\mathcal{O}$-modules $M$ equipped with a $\phi$-semilinear map $M\to
M$ such that the induced map $\phi^* M\to M$ is an isomorphism. Then there
is a functor
$$T:\Phi\operatorname{M}_{\mathcal{O}_\mathcal{E}}\to\Rep_{Z_p}(G_{K_\infty})$$
which is in fact an equivalence of abelian categories (see section 1.1.12 of
\cite{kis04}).
 Let $A$ be a finite $\Z_p$-algebra, and let
$\Rep'_A(G_{K_\infty})$ denote the category of continuous representations of
$G_{K_\infty}$ on finite $A$-modules, and let $\Rep_A(G_{K_\infty})$ denote
the full subcategory of objects which are free as $A$-modules. Let $\Phi
\operatorname{M}_{\mathcal{O}_\mathcal{E},A}$ denote the category whose
objects are objects of $\Phi \operatorname{M}_{\mathcal{O}_\mathcal{E}}$
equipped with an action of $A$.

\begin{lemma}\label{21}The functor $T$ above induces an equivalence of abelian
categories $$T_A:\Phi
\operatorname{M}_{\mathcal{O}_\mathcal{E},A}\to\Rep'_A(G_{K_\infty}).$$The
functor $T_A$ induces a functor
$$T_{\bigs,A}:(\Mod\FI/\bigs)_A\to\Rep_A(G_{K_\infty});\ \mathfrak{M}\mapsto
T_A(\mathcal{O}_\mathcal{E}\otimes_\bigs\mathfrak{M}).$$\end{lemma}
\begin{proof} Lemmas 1.2.7 and 1.2.9 of \cite{kis04}.\end{proof}

Fix $\mathbb{F}$ a finite extension of $\F_p$, and a continuous
representation of $G_K$ on a $2$-dimensional $\F$-vector space $V_{\F}$. We
suppose that $V_{\F}$ is the generic fibre of a finite flat group scheme,
and let $M_\F$ denote the preimage of $V_\F(-1)$ under the equivalence
$T_\F$ of Lemma \ref{21}.

In fact, from now on we assume that the action of $G_K$ on $V_{\F}$ is
trivial, that $k\subset\F$, and that $k\neq\F_p$. In applications we will
reduce to this case by base change.

Recall from Corollary 2.1.13 of \cite{kis04} that we have a projective
scheme $\mathcal{GR}_{V_{\F},0}$, such that for any finite extension $\F'$
of $\F$, the set of isomorphism classes of finite flat models of
$V_{\F'}=V_{\F}\otimes_{\F}\F'$ is in natural bijection with
$\mathcal{GR}_{V_{\F},0}(\F')$. We work below with the closed subscheme
$\mathcal{GR}^{\mathbf{v}}_{V_{\F},0}$ of $\mathcal{GR}_{V_{\F},0}$, defined
in Lemma 2.4.3 of \cite{kis04}, which parameterises isomorphism classes of
finite flat models of $V_{\F'}$ with cyclotomic determinant.

As in section 2.4.4 of \cite{kis04}, if $\F^{\text{sep}}$ is the residue
field of $K_0^{\text{sep}}$, and $\sigma\in\Gal(K_0/\Q_p)$, we denote by
$\epsilon_\sigma\in k\otimes_{\F_p}\F'$ the idempotent which is $1$ on the
kernel of the map $1\otimes\sigma:k\otimes_{\F_p}\F'\to\F^{\text{sep}}$
corresponding to $\sigma$, and $0$ on the other maximal ideals of
$k\otimes_{\F_p}\F'$.

\begin{lemma}\label{251}If $\F'$ is a finite extension of $\F$, then the elements of
$\mathcal{GR}^{\mathbf{v}}_{V_{\F},0}(\F')$ naturally correspond to free
$k\otimes_{\F_p}\F'[[u]]$-submodules $\mathcal{M}_{\F'}\subset
M_{\F'}:=M_\F\otimes_\F\F'$ of rank $2$ such that:\begin{enumerate}\item
$\mathcal{M}_{\F'}$ is $\phi$-stable. \item For some (so any) choice of
$k\otimes_{\F_p}\F'[[u]]$-basis for $\mathcal{M}_{\F'}$, for each $\sigma$
the map
$$\phi:\epsilon_{\sigma}\mathcal{M}_{\F'}\to\epsilon_{\sigma\circ\phi^{-1}}\mathcal{M}_{\F'}$$
has determinant $\alpha u^e$, $\alpha\in\F'[[u]]^\times$.\end{enumerate}
\end{lemma}
\begin{proof}This follows just as in the proofs of Lemma 2.5.1 and Proposition 2.2.5 of
\cite{kis04}. More precisely, the method of the proof of Proposition 2.2.5
of \cite{kis04}  allows one to ``decompose'' the determinant condition into
the condition that for each $\sigma$ we have
$$\dim_{\F'}(\epsilon_{\sigma\circ\phi^{-1}}\mathcal{M}_{\F'}/\phi(\epsilon_{\sigma}\mathcal{M}_{\F'}))=e,$$and
then an identical argument to that in the proof of Lemma 2.5.1 \cite{kis04}
shows that this condition is equivalent to the stated one.
\end{proof}
We now number the elements of $\Gal(K_0/\Q_p)$ as $\sigma_1,\dots,\sigma_r$,
in such a way that $\sigma_{i+1}=\sigma_i\circ\phi^{-1}$ (where we identify
$\sigma_{r+1}$ with $\sigma_1$). For any sublattice $\mathfrak{M}_\F$ in
$(Mod/\mathfrak{S})_\F$ and any $(A_1,\dots,A_r)\in\M_2(\F((u)))^r$, we
write $\mathfrak{M}_\F\sim A$ if there exist bases
$\{\textbf{e}^i_1,\textbf{e}^i_2\}$ for
$\epsilon_{\sigma_i}\mathcal{M}_\F$ such that $$\phi\left(%
\begin{array}{c}
 \textbf{e}^i_1 \\ \textbf{e}^i_2 \\
\end{array}%
\right)=A_i\left(%
\begin{array}{c}
 \textbf{e}^{i+1}_1  \\ \textbf{e}^{i+1}_2 \\
\end{array}%
\right).$$ If we have fixed such a choice of basis, then for any
$(B_1,\dots,B_r)\in \GL_2(k_r((u)))^r$ we denote by
$B\mathfrak{M}$ the module generated by $\left<B_i\left(%
\begin{array}{c}
 \textbf{e}^i_1 \\ \textbf{e}^i_2 \\
\end{array}%
\right)\right>$, and consider $B\mathfrak{M}$ with respect to the basis
given by these entries.

\begin{prop}\label{connected}Let $\F'/\F$ be a finite extension. Suppose
that $x_1,x_2\in\mathcal{GR}^{\mathbf{v}}_{V_\F,0}(\F')$ and that the
corresponding objects of $(\Mod/\S)_{\F'}$, $\mathfrak{M}_{\F',1}$ and
$\mathfrak{M}_{\F',2}$ are both non-ordinary. Then (the images of) $x_1$ and
$x_2$ both lie on the same connected component of
$\mathcal{GR}^{\mathbf{v}}_{V_\F,0}(\F')$.\end{prop}
\begin{proof}Replacing $V_\F$ by $\F'\otimes_\F v_\F$, we may assume that
$\F'=\F$. Suppose that $\mathfrak{M}_{\F,1}\sim A$. Then
$\mathfrak{M}_{\F,2}=B\cdot\mathfrak{M}_{\F,1}$ for some
$B\in\GL_2(k_r((u)))^r$, and $\mathfrak{M}_{\F,2}\sim(\phi(B_i)\cdot
A_i\cdot B_{i+1}^{-1})$. Each $B_i$ is uniquely determined up to left
multiplication by elements of $\GL_2(\F[[u]])$, so by the Iwasawa
decomposition we may assume that each $B_i$ is upper triangular. By Lemma
\ref{251}, $\det\phi(B_{i})\det B_{i+1}^{-1}\in\F[[u]]^\times$ for all $i$,
which implies that $\det(B_i)\in\F[[u]]^\times$ for all $i$, so that the
diagonal elements of $B_i$ are $\mu^i_1 u^{-a_i}$, $\mu^i_2 u^{a_i}$ for
$\mu^i_1$, $\mu^i_2\in\F[[u]]^\times$, $a_i\in\Z$. Replacing $B_i$ with
$\diag(\mu^i_1,\mu^i_2)^{-1}B_i$, we may assume that $B_i$ has diagonal
entries $u^{-a_i}$ and $u^{a_i}$.

We now show that $x_1$ and $x_2$ are connected by a chain of rational
curves, using the following lemma:

\begin{lemma}\label{path}Suppose that $(N_i)$ are nilpotent elements of $M_2(\F((u)))$
such that $\mathfrak{M}_{\F,2}=(1+N)\cdot\mathfrak{M}_{\F,1}$. If
$\phi(N_{i})AN_{i+1}^{\ad}\in M_2(F[[u]])$ for all $i$, then there is a map
$\mathbb{P}^1\to\mathcal{GR}^{\mathbf{v}}_{V_\F,0}$ sending $0$ to $x_1$ and
$1$ to $x_2$.
\end{lemma}
\begin{proof}Exactly as in the proof of Lemma 2.5.7 of \cite{kis04}.
\end{proof}In fact, we will only apply this lemma in situations where all
but one of the $N_i$ are zero, so that the condition of the lemma is
automatically satisfied.

\begin{lemma}With respect to some basis, $\phi:M_\F\to M_\F$ is given by
$\bigl(\begin{smallmatrix}1&0\\0&1\end{smallmatrix}\bigr)$.\end{lemma}
\begin{proof}This is immediate from the definition of $M_\F$ (recall that we
have assumed that the action of $G_K$ on $V_\F$ is trivial).
\end{proof}
Let $v_1$, $v_2$ be a basis as in the lemma, and let $\mathfrak{M}_\F$ be
the sub-$k\otimes_{\F_p}\F[[u]]$-module generated by $u^{e/(p-1)}v_1$ and
$v_2$ (note that the assumption that the action of $G_K$ on $V_\F$ is
trivial guarantees that $e|(p-1)$). Then $\mathfrak{M}_\F$ corresponds to an
object of $\mathcal{GR}^{\mathbf{v}}_{V_\F,0}(\F')$, and
$\mathfrak{M}_\F\sim(A_i)$ where each
$A_i=\bigl(\begin{smallmatrix}u^e&0\\0&1\end{smallmatrix}\bigr)$, so that
$\mathfrak{M}_\F$ is ordinary.

Furthermore, every object of $\mathcal{GR}^{\mathbf{v}}_{V_\F,0}(\F')$ is
given by $B\cdot\mathfrak{M}_\F$ for some $B=(B_i)$, where
$B_i=\bigl(\begin{smallmatrix}u^{-a_i}&v_i\\0&u^{a_i}\end{smallmatrix}\bigr)$,
and $\phi(B_i)A_i B_{i+1}^{-1}\in M_2(\F[[u]])$ for all $i$. Examining the
diagonal entries of $\phi(B_i) A_i B_{i+1}^{-1}$, we see that we must have
$e\geq pa_{i}-a_{i+1}\geq 0$ for all $i$.
\begin{lemma}We have $e/(p-1)\geq a_i\geq 0$ for all $i$. Furthermore, if
any $a_i=0$ then all $a_i=0$; and if any $a_i=e/(p-1)$, then all
$a_i=e/(p-1)$.\end{lemma}
\begin{proof}Suppose that $a_j\leq 0$. Then $pa_j\geq a_{j+1}$, so
$a_{j+1}\leq 0$. Thus $a_i\leq 0$ for all $i$. But adding the inequalities
gives $(p-1)(a_1+\dots+a_r)\geq 0$, so in fact $a_1=\dots=a_r=0$. The other
half of the lemma follows in a similar fashion.
\end{proof}
Note that the ordinary objects are precisely those with all $a_i=0$ or all
$a_i=e/(p-1)$. We now show that there is a chain of rational curves linking
any non-ordinary point to the point corresponding to
$C\cdot\mathfrak{M}_\F$, where
$C_i=\bigl(\begin{smallmatrix}u^{-1}&0\\0&u\end{smallmatrix}\bigr)$.

Choose a non-ordinary point $D\cdot\mathfrak{M}_\F$,
$D_i=\bigl(\begin{smallmatrix}u^{-b_i}&w_i\\0&u^{b_i}\end{smallmatrix}\bigr)$.
We claim that there is a chain of rational curves linking this to the point
$D'\cdot\mathfrak{M}_\F$,
$D'_i=\bigl(\begin{smallmatrix}u^{-b_i}&0\\0&u^{b_i}\end{smallmatrix}\bigr)$.
Clearly, it suffices to demonstrate that there is a rational curve from
$D\cdot\mathfrak{M}_\F$ to the point $D^j\cdot\mathfrak{M}_\F$, where
$$D^j_i=\left\{\begin{array}{ll}D_i,& i\neq
j\\ \left(\begin{smallmatrix}u^{-b_j}&0\\0&u^{b_j}\end{smallmatrix}\right),&
i=j.\end{array}\right. $$ But this is easy; just apply Lemma \ref{path} with
$N=(N_i)$, $$N_i=\left\{\begin{array}{ll}0,& i\neq j\\
\left(\begin{smallmatrix}0&-w_ju^{-b_j}\\0&0\end{smallmatrix}\right),&
i=j.\end{array}\right. $$It now suffices to show that there is a chain of
rational curves linking $D'\cdot\mathfrak{M}_\F$ to $C\cdot\mathfrak{M}_\F$.
Suppose that  $D''\cdot\mathfrak{M}_\F$ also corresponds to a point of
$\mathcal{GR}^{\mathbf{v}}_{V_\F,0}(\F)$, where for some $j$ we have
$$D''_i=\left\{\begin{array}{ll}D'_i,& i\neq j\\
\left(\begin{smallmatrix}u^{1-b_j}&0\\0&u^{b_j-1}\end{smallmatrix}\right),&
i=j.\end{array}\right. $$ Then we claim that there is a rational curve
linking $D'\cdot\mathfrak{M}_\F$ and $D''\cdot\mathfrak{M}_\F$. Note that
$D''=ED'$, where $$E_i=\left\{\begin{array}{ll}1,& i\neq j\\
\left(\begin{smallmatrix}u^{-1}&0\\0&u\end{smallmatrix}\right),&
i=j.\end{array}\right. $$ Since
$\left(\begin{smallmatrix}u^{-1}&0\\0&u\end{smallmatrix}\right)=\left(\begin{smallmatrix}0&1\\-1&2u\end{smallmatrix}\right)\left(\begin{smallmatrix}2&-u\\u^{-1}&0\end{smallmatrix}\right)$,
and
$\left(\begin{smallmatrix}0&1\\-1&2u\end{smallmatrix}\right)\in\GL_2(\F[[u]])$,
we can apply Lemma \ref{path} with $$N_i=\left\{\begin{array}{ll}0,& i\neq j\\
\left(\begin{smallmatrix}1&-u\\u^{-1}&-1\end{smallmatrix}\right),&
i=j.\end{array}\right. $$ Proposition \ref{connected} now follows from:
\begin{lemma}If $e/(p-1)>a_i>0$ and $e\geq pa_{i}-a_{i+1}\geq 0$ for all $i$,
and not all the $a_i$ are equal to $1$, then for some $j$ we can define
$$a'_i=\left\{\begin{array}{ll}a_i,& i\neq
j\\a_j-1,& i=j\end{array}\right.$$and we have $e\geq pa'_{i}-a'_{i+1}\geq 0$
for all $i$.\end{lemma}
\begin{proof}Suppose not. Then for each $i$, either $pa_{i-1}-(a_{i}-1)>e$, or
$p(a_i-1)-a_{i+1}<0$; that is, either $pa_{i-1}-a_{i}=e$, or $p-1\geq
pa_i-a_{i+1}\geq 0$. Comparing the statements for $i$, $i+1$, we see that
either $pa_{i}-a_{i+1}=e$ for all $i$, or $p-1\geq pa_i-a_{i+1}\geq 0$ for
all $i$. In the former case we have $a_i=e/(p-1)$ for all $i$, a
contradiction. In the latter case, summing the inequalities gives
$r(p-1)\geq(p-1)(a_1+\dots+a_r)\geq (r+1)(p-1)$, a
contradiction.\end{proof}\end{proof}

\section{Modularity lifting theorems}The results of section \ref{ctd} can
easily be combined with the machinery of \cite{kis04} to yield modularity
lifting theorems. For example, we have the following:

\begin{thm}\label{41}Let $p>2$, let $F$ be a totally real field, and
let $E$ be a finite extension of $\Q_p$ with ring of integers $\mathcal{O}$.
Let $\rho:G_F\to\GL_2(\mathcal{O})$ be a continuous representation
unramified outside of a finite set of primes, with determinant a finite
order character times the $p$-adic cyclotomic character. Suppose that
\begin{enumerate}\item $\rho$ is potentially Barsotti-Tate at each $v|p$.
\item There exists a Hilbert modular form $f$ of parallel weight 2 over $F$ such that $\overline{\rho}_f\sim\overline{\rho}$, and for each $v|p$, $\rho$ is potentially ordinary at $v$ if and only if $\rho_f$ is. \item $\overline{\rho}|_{F(\zeta_p)}$ is
absolutely irreducible, and if $p>3$ then
$[F(\zeta_p):F]>3$.\end{enumerate}Then $\rho$ is modular.\end{thm}
\begin{proof}The proof of this theorem is almost identical to the proof of
Theorem 3.5.5 of \cite{kis04}. Indeed, the only changes needed are to
replace property \emph{(iii)} of the field $F'$ chosen there by
``\emph{(iii)} If $v|p$ then $\overline{\rho}|G_{F_v}$ is trivial, and the
residue field at $v$ is not $\F_p$'', and to note that Theorem 3.4.11 of
\cite{kis04} is still valid in the context in which we need it, by
Proposition \ref{connected}.\end{proof}

For the applications to mod $p$ Hilbert modular forms in \cite{gee053} it is
important not to have to assume that $\rho$ is potentially ordinary at $v$
if and only if $\rho_f$ is. Fortunately, in \cite{gee053} it is only
necessary to work with totally real fields in which $p$ is unramified, and
in that case we are able, following \cite{kis05}, to remove this assumption.

\begin{thm}Let $p>2$, let $F$ be a totally real field in which $p$ is unramified, and
let $E$ be a finite extension of $\Q_p$ with ring of integers $\mathcal{O}$.
Let $\rho:G_F\to\GL_2(\mathcal{O})$ be a continuous representation
unramified outside of a finite set of primes, with determinant a finite
order character times the $p$-adic cyclotomic character. Suppose that
\begin{enumerate}\item $\rho$ is potentially Barsotti-Tate at each $v|p$.
\item $\overline{\rho}$ is modular. \item $\overline{\rho}|_{F(\zeta_p)}$ is
absolutely irreducible.\end{enumerate}Then $\rho$ is modular.\end{thm}
\begin{proof}Firstly, note that by a standard result (see e.g. \cite{bdj}) we have $\overline{\rho}\sim\overline{\rho}_f$, where $f$ is a form of parallel
weight $2$. We now follow the proof of Corollary 2.13 of \cite{kis05}. Let
$\mathcal{S}'$ denote the set of $v|p$ such that $\rho|_{G_v}$ is
potentially ordinary. After applying Lemma
\ref{existord} below, we may assume that
$\overline{\rho}\sim\overline{\rho}_f$, where $f$ is a form of parallel
weight $2$, and $\rho_f$ is potentially ordinary and potentially
Barsotti-Tate for all $v\in\mathcal{S}'$.

We may now make a solvable base change so that the hypotheses on $F$ in
Theorem \ref{41} are still satisfied, and in addition $[F:\Q]$ is even, and
at every place $v|p$ $f$ is either unramified or special of conductor $1$. By our choice of $f$, $\rho_f|_{G_v}$ is Barsotti-Tate and ordinary at each place $v\in \mathcal{S}'$. In order to apply Theorem \ref{41}, we need to check that we can replace $f$ by a form $f'$ such that
$\overline{\rho}\sim\overline{\rho}_{f'}$, and $\rho_{f'}$ is Barsotti-Tate
at all $v|p$ and is ordinary if and only if $\rho$ is. That is, we wish to choose $f'$ so that $\rho_{f'}$ is Barsotti-Tate
and ordinary at all places $v\in\mathcal{S}'$, and is Barsotti-Tate
and non-ordinary at all other places dividing $p$. 
The existence of such an $f'$ follows at once from the proof of Theorem 3.5.7 of \cite{kis04}. The theorem then
follows from Theorem \ref{41}.
\end{proof}
\begin{lemma}\label{existord}Let $F$ be a totally real field in which $p$ is
unramified, and $\mathcal{S}'$ a set of places of $F$ dividing $p$. Let $f$
be a Hilbert modular cusp form over $F$ of parallel weight $2$, with
$\overline{\rho}_f$ absolutely irreducible, and suppose that for
$v\in\mathcal{S}'$ $\overline{\rho}_f|_{G_{F_v}}$ is the reduction of a
potentially Barsotti-Tate representation of $G_{F_v}$ which is also
potentially ordinary.

Then there is a Hilbert modular cusp form $f'$ over $F$ of parallel weight
$2$ with $\overline{\rho}_{f'}\sim\overline{\rho}_f$, and such that for all
$v\in\mathcal{S}'$, $\rho_{f'}$ becomes ordinary and Barsotti-Tate over some
finite extension of $F_v$.\end{lemma}
\begin{proof}We follow the proof of Lemma 2.14 of \cite{kis05}, indicating
only the modifications that need to be made. Replacing the appeal to
\cite{cdt} with one to Proposition 1.1 of \cite{dia05}, the proof of Lemma 2.14 of 
\cite{kis05} shows that we can find $f'$ such that
$\overline{\rho}_{f'}\sim\overline{\rho}_f$, and such that for all
$v\in\mathcal{S}'$, $\rho_{f'}$ becomes Barsotti-Tate over
$F_v(\zeta_{q_v})$, where $q_v$ is the degree of the residue field of $F$ at
$v$. Furthermore, we can assume that the type of $\rho_{f'}|_{G_{F_v}}$ is
$\widetilde{\omega}_1\oplus\widetilde{\omega}_2$, where
$\overline{\rho}_{f'}|_{G_{F_v}}\sim\left(\begin{smallmatrix}\omega_1\chi&*\\0&\omega_2\end{smallmatrix}\right)$,
where $\chi$ is the cyclotomic character, and a tilde denotes the
Teichmuller lift. Let $\mathcal{G}$ denote the $p$-divisible group over
$\mathcal{O}_{F_v(\zeta_{q_v})}$ corresponding to
$\rho_{f'}|_{F_v(\zeta_{q_v})}$ Then by a scheme-theoretic closure argument,
$\mathcal{G}[p]$ fits into a short exact sequence
$$0\to\mathcal{G}_1\to\mathcal{G}[p]\to\mathcal{G}_2\to 0.$$The information
about the type then determines the descent data on the Breuil modules
corresponding to $\mathcal{G}_1$ and $\mathcal{G}_2$. We will be done if we
can show that $\mathcal{G}_1$ is multiplicative and $\mathcal{G}_2$ is
\'{e}tale. However, by the hypothesis on $\mathcal{S}'$  we can write down a multiplicative group
scheme $\mathcal{G}'_1$ with the same descent data and generic fibre as
$\mathcal{G}_1$. Then Lemma \ref{final} below shows that $\mathcal{G}_1$ is indeed multiplicative.
The same argument shows that $\mathcal{G}_2$ is  \'{e}tale.
\end{proof}

\begin{lemma}\label{final}Let $k$ be a finite field of characteristic $p$, and let $L=W(k)[1/p]$. Fix $\pi=(-p)^{1/(p^d-1)}$ where $d=[k:\F_p]$, and let $K=L(\pi)$. Let $E$ be a finite field containing $k$. Let $\mathcal{G}$ and $\mathcal{G}'$ be finite flat rank one $E$-module schemes over $\mathcal{O}_K$ with generic fibre descent data to $L$. Suppose that the generic fibres of  $\mathcal{G}$ and $\mathcal{G}'$ are isomorphic as $G_L$-representations, and that $\mathcal{G}$ and $\mathcal{G}'$ have the same descent data. Then $\mathcal{G}$ and $\mathcal{G}'$ are isomorphic.\end{lemma}
\begin{proof}This follows from a direct computation using Breuil modules with descent data. Specifically, it follows at once from Example A.3.3 of \cite{sav06}, which computes the generic fibre of any  finite flat rank one $E$-module scheme over $\mathcal{O}_K$ with generic fibre descent data to $L$. \end{proof}

	\section{Erratum}
	Unfortunately, the proof of Theorem 3.2  (and thus the main Theorem stated in the introduction to the paper) is incomplete; in particular, the proof of Lemma 3.3 is incorrect. Specifically, one cannot automatically assume that the type of $\rho_{f'}$ is $\tilde{\omega}_{1}\oplus\tilde{\omega}_{2}$; to do so is to make a rather strong assumption about the Serre weights of $\rhobar_{f}$. In addition one cannot conclude that the type determines the descent data on $\mathcal{G}_{1}$ and $\mathcal{G}_{2}$, at least in the case where $\rhobar_{f}|_{G_{F_{v}}}$ is split; either $\mathcal{G}_{1}$ can correspond to $\tilde{\omega}_{1}$ and $\mathcal{G}_{2}$ to $\tilde{\omega}_{2}$, or vice versa.

	As a consequence, we are only able to obtain a slightly weaker modularity lifting theorem; the most general result we can obtain is:

	\begin{thm}Let $p>2$, let $F$ be a totally real field, and
	let $E$ be a finite extension of $\Q_p$ with ring of integers $\mathcal{O}$.
	Let $\rho:G_F\to\GL_2(\mathcal{O})$ be a continuous representation
	unramified outside of a finite set of primes, with determinant a finite
	order character times the $p$-adic cyclotomic character. Suppose that
	\begin{enumerate}\item $\rho$ is potentially Barsotti-Tate at each $v|p$.
	\item There exists a Hilbert modular form $f$ of parallel weight 2 over $F$ such that $\overline{\rho}_f\sim\overline{\rho}$, and for each $v|p$, if $\rho$ is potentially ordinary at $v$ then so is $\rho_f$. \item $\overline{\rho}|_{F(\zeta_p)}$ is
	absolutely irreducible, and if $p=5$ and the projective image of $\rhobar$ is isomorphic to $\operatorname{PGL}_{2}(\F_{5})$ then
	$[F(\zeta_5):F]=4$.\end{enumerate}Then $\rho$ is modular.\end{thm}
	\begin{proof}The proof is extremely similar to that of Theorem 3.1. Hypothesis (3) has been weakened because of a corresponding weakening of (3.2.3)(3) in the final version of \cite{kis04}. As for (2), we need only check that after making a base change, we may assume that at each place dividing $p$, $\rho$ is potentially ordinary if and only if $\rho_f$ is. This is easily achieved by employing Lemma 3.1.5 of \cite{kis04} at each place where $\rho$ is not potentially ordinary.\end{proof}

	Additionally, we would like to thank Fred Diamond and Florian Herzig for independently bringing  to our attention a minor error in the proof of Proposition 2.3. The points $D^j$ constructed in the proof are not necessarily points on $\mathcal{GR}_{V_{\F},0}$. However, their only use is in showing that the points $D$ and $D'$ lie on the same component, and this in fact follows immediately from an application of Lemma 2.4 with $N=(N_i)$, $$N_i=
	\left(\begin{smallmatrix}0&-w_iu^{-b_i}\\0&0\end{smallmatrix}\right).$$

\def\cprime{$'$}
\providecommand{\bysame}{\leavevmode\hbox to3em{\hrulefill}\thinspace}
\providecommand{\MR}{\relax\ifhmode\unskip\space\fi MR }
\providecommand{\MRhref}[2]{%
  \href{http://www.ams.org/mathscinet-getitem?mr=#1}{#2}
}
\providecommand{\href}[2]{#2}

\end{document}